\documentclass[a4paper,12pt]{article}
\usepackage{amsmath}
\usepackage[mathscr]{eucal}
\usepackage{amssymb}
\usepackage{latexsym}
\usepackage{amsthm}
\usepackage{amscd}
\theoremstyle{plain}
\newtheorem{theorem}{Theorem}

\newtheorem{lemma}{Lemma}[section]

\theoremstyle{remark}
\newtheorem{remark}{Remark}
\theoremstyle{definition}

\usepackage[dvips]{graphicx}
\usepackage{psfrag}
\begin{document}

\title{ Infinitesimal Torelli theorem for surfaces with \\ 
$c_1^2=3$, $\chi = 2$, 
and the torsion group $\mathbb{Z}/3$}
\author{Masaaki Murakami \footnote{The author was supported 
by INDAM-GNSAGA during his stay at University 
of Padua, Italy, starting from October 2003 to March 2004, 
in which he performed all the 
computations included in this article.}}
\date{}
\maketitle
\begin{abstract} 
We prove the infinitesimal Torelli theorem
for general minimal complex surfaces $X$'s  
with the first Chern number $3$, 
the geometric genus $1$, and the irregularity $0$  
which have non-trivial $3$-torsion divisors. 
We also show that the coarse moduli space for 
surfaces with the invariants as above is 
a $14$-dimensional unirational variety. 
\end{abstract}
\setcounter{section}{-1}
\section{Introduction}  \label{sectn:introduction}

In the present paper, we will prove the infinitesimal 
Torelli theorem for general minimal complex surfaces $X$'s 
with $c_1^2 = 3$, $\chi (\mathcal{O}) =2$, and 
$\mathrm{Tors} (X) \simeq \mathbb{Z} /3$, where 
$c_1$, $\chi (\mathcal{O})$, and $\mathrm{Tors} (X)$ are 
the first Chern class, the Euler characteristic of the 
structure sheaf, and the torsion part of the Picard group
of $X$, respectively. 
We will also show that all surfaces with the invariants as above 
are deformation equivalent to each other, 
and that their coarse moduli space  $\mathcal{M}$  
is a $14$-dimensional unirational variety. 
Here, let us remark that 
the condition $\mathrm{Tors} (X) \simeq \mathbb{Z} /3$ is 
a topological one; minimal surfaces with 
$c_1^2 = 3$ and $\chi (\mathcal{O}) =2$ have  
the geometric genus $p_g=1$ and the irregularity $q=0$, 
hence the torsion group $\mathrm{Tors} (X)$ isomorphic 
to the first homology group $H_1(X, \mathbb{Z})$.  

As is well known today, Torelli type theorems do not necessarily 
hold for surfaces. One of the most famous counter examples  
are surfaces of general type with $p_g = q =0$. 
Although Torelli type theorems have been proved for many classes 
of surfaces, finding what conditions we should impose still remains 
as a challenging problem. So it makes sense to study period maps for 
concrete classes of surfaces. 
 
Let us recall some results on period maps for surfaces of general 
type with $p_g=1$ and $q=0$. 
In \cite{pg1k2=1}, Catanese proved the infinitesimal Torelli 
theorem for general minimal surfaces with $c_1^2 =1$, $p_g=1$,
and $q=0$, while in \cite{globalpg1k2=1} that the global period mapping has 
degree at least $2$. He first showed that any such surface is essentially 
a weighted complete intersection of type $(6,6)$ in the weighted 
projective space $\mathbb{P} (1,2,2,3,3)$, and used this complete 
description to study the period map for these surfaces. 
Meanwhile for the case $c_1^2 =2$, $p_g=1$, and $q=0$, 
the torsion group is either $0$ or $\mathbb{Z} /2$.  
Using a complete description for the case 
of $\mathbb{Z} /2$    
by Catanese and Debarre \cite{pg1k2=2}, 
Oliverio studied in \cite{torellik2=2} 
the infinitesimal period maps for the case of 
non-trivial $2$-torsion divisors by the same method as in 
\cite{pg1k2=1}. 

Consider the case $c_1^2 =3$. 
In this case, the order $\sharp \mathrm{Tors} (X)$ is at most 
$3$ by a result in \cite{bound}. 
Moreover, in \cite{3tors'}, the author showed that 
any surface $X$ of this class with  
$\mathrm{Tors} (X) \simeq \mathbb{Z} /3$ is essentially 
a quotient of a $(3,3)$-complete intersection 
in $\mathbb{P}^4$ by a certain free action by $\mathbb{Z} /3$. 
Using this complete description, we will show 
in the present paper the infinitesimal 
Torelli theorem for general $X$'s by the same method as in 
\cite{pg1k2=1} and \cite{torellik2=2}. Here, general $X$ means any surface 
corresponding to a point in a certain Zariski open subset of 
the coarse moduli space $\mathcal{M}$.

In Section \ref{scn:results}, 
we state our main theorems of the present paper and 
recall our previous results given in \cite{3tors'}. 
In Section \ref{scn:unirat}, 
we show the unirationality of the coarse moduli space $\mathcal{M}$. 
Finally in Section \ref{scn:inftorelli}, we prove the infinitesimal 
Torelli theorem for our surfaces $X$'s.  
Throughout this paper, we work over the complex number field $\mathbb{C}$.

\medskip
{\sc Notation}
\medskip

Let $S$ be a compact complex manifold of dimension $2$. 
We denote by $p_g(S)$, $q(S)$, and $K_S$, the geometric genus, 
the irregularity and a canonical divisor of $S$, respectively.  
We denote by $\mathrm{Tors(S)}$ the torsion part of the Picard 
group, and call it the torsion group of $S$. 
For a coherent sheaf $\mathcal{F}$ on $S$, 
we denote by  $h^i(\mathcal{F})$ the dimension of 
the $i$-th cohomology group $H^i(S,\mathcal{F})$. 
The sheaf $\mathcal{O}_S$, $\varOmega^p_S$, and $\varTheta_S$ are 
the structure sheaf, the sheaf of germs of holomorphic $p$-forms, 
and that of germs of holomorphic vector fields on $S$, respectively. 
As usual, $\mathbb{P}^n$ is the projective space of dimension $n$.
We denote by $\varepsilon = \mathrm{exp}(2\pi \sqrt{-1}/3)$ 
a third root of unity.

\section{Statement of main results}  \label{scn:results}

In a previous paper \cite{3tors'}, the author gave a 
complete description for minimal algebraic surfaces $X$'s 
with $c_1^2 = 3$, $\chi (\mathcal{O}) = 2$,  
and $\mathrm{Tors} (X) \simeq \mathbb{Z} / 3$, 
where $c_1$, $\chi (\mathcal{O})$, and $\mathrm{Tors} (X)$ are 
the first Chern class, 
the Euler characteristic of the structure sheaf,   
and the torsion part of the Picard group of $X$, respectively. 
In the present paper, we will give proofs for the 
following two theorems: 

\begin{theorem} \label{thm:moduli}
All minimal algebraic surfaces $X$'s 
with $c_1^2 = 3$, $\chi (\mathcal{O}) = 2$,  
and $\mathrm{Tors} (X) \simeq \mathbb{Z} / 3$ 
are deformation equivalent to each other. 
Their coarse moduli space $\mathcal{M}$ 
is a $14$-dimensional unirational variety.   
\end{theorem} 

\begin{theorem}  \label{thm:infperiodmap}
Let $X$ be any general surface as in Theorem \ref{thm:moduli}.  
Then the infinitesimal period map 
$\mu : H^1 (\varTheta_X) \to 
\mathrm{Hom} (H^0(\varOmega_X^2), H^1(\varOmega_X^1))$ 
is injective. 
\end{theorem}

\begin{remark}
The surfaces $X$'s as in Theorem \ref{thm:moduli}
have the geometric genus $p_g=1$ and the irregularity $q=0$. 
We refer the readers to \cite{globalmoduli} 
for the existence of the coarse moduli space $\mathcal{M}$. 
See also \cite{periodmapgriff} for the infinitesimal period map.
\end{remark}

In order to give proofs for the theorems above, 
let us first recall the main results given in \cite{3tors'}. 
See \cite{3tors'} for proofs of the following two theorems:

\begin{theorem}[\cite{3tors'}] \label{thm:description}
Let $X$ be a minimal algebraic surface with 
$c_1^2 = 3$, $\chi (\mathcal{O}) =2$, 
and $\mathbb{Z} / 3 \subset \mathrm{Tors} (X)$. 
Let $\pi : Y \to X$ be the unramified Galois triple cover 
corresponding to a non-trivial $3$-torsion divisor. 
Then both the fundamental group $\pi_1 (X)$ and the 
torsion group $\mathrm{Tors} (X)$ are isomorphic to 
the cyclic group $\mathbb{Z} /3$. Further, the canonical 
model $Z$ of $Y$ is a complete intersection in the 
$4$-dimensional projective space $\mathbb{P}^4$ defined 
by two homogeneous polynomials $\tilde{F_1}$ and $\tilde{F_2}$ 
of degree $3$ satisfying 
\[
  \tilde{F_i} 
  (W_0, \varepsilon X_1, \varepsilon X_2, 
        \varepsilon^{-1} Y_3, \varepsilon^{-1} Y_4) 
  = \tilde{F_i}
    (W_0, X_1, X_2, Y_3, Y_4) \quad (i=1,2). 
\]
Here, $(W_0, X_1, \cdots, Y_4)$ is 
a homogeneous coordinate of $\mathbb{P}^4$,  
and the constant $\varepsilon = \exp (2\pi \sqrt{-1} /3)$ 
is a third root of unity.  
\end{theorem}

\begin{theorem}[\cite{3tors'}] \label{thm:kuranishisp}
  Let $X$ be a surface as in Theorem \ref{thm:description}. 
If $X$ has an ample canonical divisor $K_X$, 
then $h^1(\varTheta_X) = 14$ and $h^2(\varTheta_X) = 0$,  
hence the Kuranishi space of $X$ is smooth and of dimension $14$. 
\end{theorem}

\begin{remark} \label{rem:expliciteforms}
Explicit forms of the two polynomials in 
Theorem \ref{thm:description} are given by 
\begin{equation} \label{eqn:expliciteforms}
 \tilde{F_i} = a_0^{(i)} W_0^3 
              + W_0 \tilde{\gamma_i} (X_1, X_2, Y_3, Y_4)
              + \tilde{\alpha_i} (X_1, X_2) 
              + \tilde{\beta_i}  (Y_3, Y_4) 
\end{equation}
for $i = 1,2$, where 
\begin{align*}
 \tilde{\gamma_i} &= a_1^{(i)} X_1 Y_3  +  a_2^{(i)} X_1 Y_4 + 
                     a_3^{(i)} X_2 Y_3  +  a_4^{(i)} X_2 Y_4 ,  \\
 \tilde{\alpha_i} &= a_5^{(i)} X_1^3    +  a_6^{(i)} X_1^2 X_2 +
                     a_7^{(i)} X_1 X_2^2  +  a_8^{(i)} X_2^3 ,  \\
 \tilde{\beta_i}  &= a_9^{(i)} Y_3^3      +  a_{10}^{(i)} Y_3^2 Y_4 +
                     a_{11}^{(i)} Y_3 Y_4^2  + a_{12}^{(i)} Y_4^3 ,  
\end{align*}
are homogeneous polynomials of $X_1, \cdots , Y_4$
with coefficients $a_j^{(i)} \in \mathbb{C}$. 
\end{remark}

\begin{remark} \label{rem:galoisgp}
 The complete intersection $Z$ is the image of the 
canonical map $\varPhi_{K_Y} : Y \to \mathbb{P}^4$.  
We have a natural action on $Z$ by 
the Gaolis group $G=\mathrm{Gal} (Y/X) \simeq \mathbb{Z} /3$ 
of $Y$ over $X$. 
This action is given by 
\begin{equation} \label{eqn:action}
  \tau_0 : (W_0 : X_1 : X_2 : Y_3 : Y_4) 
             \mapsto 
           (W_0 : \varepsilon X_1 : \varepsilon X_2 : 
                 \varepsilon^{-1} Y_3 : \varepsilon^{-1} Y_4) , 
\end{equation}
where $\tau_0$ is a generator of the group $G$. 
Since this action on $Z$ has no fixed points, 
the coefficients $a_j^{(i)}$'s satisfy the 
following three conditions:

i) at least one out of $a_0^{(1)}$ and $a_0^{(2)}$ are not equal to zero,

ii) two polynomials $\tilde{\alpha_1}$ and $\tilde{\alpha_2}$ 
    have no common zeroes on $\mathbb{P}^1= \{ (X_1 : X_2) \}$, 

iii) two polynomials $\tilde{\beta_1}$ and $\tilde{\beta_2}$ 
    have no common zeroes on $\mathbb{P}^1= \{ (Y_3 : Y_4) \}$.  

\end{remark}

For each integer $n \geq 0$, 
we have a natural isomorphism 
\begin{equation} \label{eql:decomposition}
 H^0(\mathcal{O}_Y (nK_Y)) \simeq 
  \bigoplus_{m = 0, 1, -1} H^0(\mathcal{O}_X (nK_X - mT_0))
\end{equation}
corresponding to the action by $G$, 
where $T_0$ is a generator of the torsion group $\mathrm{Tors} (X)$. 
Note that this is a decomposition into homogeneous eigen spaces, 
and that, in Theorem \ref{thm:description}, 
the sets $\{ W_0 \}$, $\{ X_1, X_2\}$, and $\{ Y_3, Y_4 \}$ 
correspond to a base of $H^0(\mathcal{O}_X (K_X))$, 
that of $H^0(\mathcal{O}_X (K_X - T_0))$, 
and that of $H^0(\mathcal{O}_X (K_X + T_0))$, respectively. 
The polynomials $\tilde{F_1}$ and $\tilde{F_2}$ generate 
the linear space consisting of all the elements in 
$H^0(\mathcal{O}_{\mathbb{P}^4}(3H))$ vanishing along $Z$, 
where $H$ is a hyperplane in $\mathbb{P}^4$. 

\section{Unirationality of the moduli space} \label{scn:unirat}

In this section, we will give a proof for 
Theorem \ref{thm:moduli}. 
We denote by $W=\mathbb{P}^4$ and $(W_0: X_1: X_2: Y_3: Y_4)$,  
the $4$-dimensional complex projective space 
and its homogeneous coordinate, respectively.   

Let $\tilde{B}$ be the set of all 
$(a_j^{(i)})_{0 \leq j \leq 12}^{1 \leq i \leq 2} \in \mathbb{C}^{26}$ 
satisfying the conditions i), ii) and iii) in Remark 
\ref{rem:galoisgp} 
such that two polynomials $\tilde{F_1}$ and $\tilde{F_2}$ given by  
(\ref{eqn:expliciteforms})  
define in $W=\mathbb{P}^4$ a complete intersection 
with at most rational double points as its singularities. 
We denote by $\tilde{B_0}$ the set of points in $\tilde{B}$ 
corresponding to non-singular complete intersections. 
Note by \cite[Remark 1]{3tors'}, we have $\tilde{B_0} \neq \emptyset$, 
hence the spaces $\tilde{B}$ and $\tilde{B_0}$ are 
dense Zariski open subsets of $\mathbb{C}^{26}$. 
We have a flat family $\tilde{\mathcal{Y}} \to \tilde{B}$ 
whose fiber on each $(a_j^{(i)}) \in \tilde{B}$ is 
a complete intersection defined by 
$\tilde{F_1}$ and $\tilde{F_2}$ 
with $a_j^{(i)}$'s as their coefficients. 
This $\tilde{\mathcal{Y}}$ is a subvariety of $\tilde{B} \times W$  
stable under the action by 
$G \simeq \langle \mathrm{id}_{\tilde{B}} \times \tau_0 \rangle 
\simeq \mathbb{Z} / 3$, where $\tau_0$ is an automorphism of 
$W$ given by (\ref{eqn:action}).
Taking a quotient of $\tilde{\mathcal{Y}}$ by this action,
we obtain a family $\tilde{\mathcal{X}} \to \tilde{B}$ 
whose fibers are the canonical models of surfaces $X$'s   
as in Theorem \ref{thm:description}. 
Note that both restrictions 
$\tilde{\mathcal{Y}} |_{\tilde{B_0}}  \to \tilde{B_0}$  
and $\tilde{\mathcal{X}} |_{\tilde{B_0}}  \to \tilde{B_0}$  
are analytic families. 

\begin{lemma}  \label{lm:coefreduction}
Let $X$ be an algebraic surface as in Theorem \ref{thm:description}. 
Then there exist bases of $H^0(\mathcal{O}_X (K_X))$, 
$H^0(\mathcal{O}_X (K_X-T_0))$, and $H^0(\mathcal{O}_X (K_X+T_0))$ 
such that the polynomials $\tilde{F_1}$ and 
$\tilde{F_2}$ satisfy $a_0^{(1)}=1$,
$a_5^{(1)}=a_9^{(1)}=1$,        $a_7^{(1)}=a_8^{(1)}=a_{12}^{(1)}=0$,
$a_8^{(2)}=a_{12}^{(2)}=1$, and $a_5^{(2)}=a_6^{(2)}=a_9^{(2)}=0$.
\end{lemma}

Proof. Take those bases and 
$\tilde{F_i}$'s in Theorem \ref{thm:description} in such a way that 
each $\tilde{\alpha_i}$ for $i=1, 2$ has 
a zero of order at least $2$ at $(X_1: X_2) = ( i-1 : 2-i )$  
and that each $\tilde{\beta_i}$ for $i=1,2$ has 
a zero at $(Y_3 : Y_4) = ( i-1 : 2-i )$.   
This is possible, since   
$(X_1 : X_2) \mapsto 
(\tilde{\alpha_1}(X_1, X_2) : \tilde{\alpha_2}(X_1, X_2))$ is 
a morphism of degree $3$. 
Then, by the conditions ii) and iii) in Remark \ref{rem:galoisgp},
we have $a_5^{(1)} \neq 0$, $a_8^{(2)} \neq 0$, 
$a_9^{(1)} \neq 0$ and $a_{12}^{(2)} \neq 0$. 
Now, by replacing the elements in these bases by 
their multiples by non-zero constants, 
and changing indices if necessary, 
we easily obtain the assertion.     \qed

Consider the case of $X$ for which $\tilde{F_i}$'s 
as in Lemma \ref{lm:coefreduction} satisfy $a_0^{(2)} \neq 0$. 
In this case, we replace $X_2$ and $Y_4$ 
by their multiples by a non-zero constant such that 
the equality $a_0^{(2)}=1$, 
as much as the equalities in the lemma above, holds. 
Then the defining polynomials $F_i = \tilde{F_{i}}$'s of 
$Z$ in $W$ are given by 
\begin{equation} \label{eqn:eqnunivcov}
 F_i =  W_0^3 
      + W_0 \gamma_i (X_1, X_2, Y_3, Y_4)
      + \alpha_i (X_1, X_2) 
      + \beta_i  (Y_3, Y_4) ,   
\end{equation}
for $i=1,2$, where  
\begin{align*}
   \gamma_1 &=  a^{(1)} X_1 Y_3 + b^{(1)} X_1 Y_4 
              + c^{(1)} X_2 Y_3 + d^{(1)} X_2 Y_4  , \\
   \gamma_2 &=  a^{(2)} X_1 Y_3 + b^{(2)} X_1 Y_4 
             + c^{(2)} X_2 Y_3  + d^{(2)} X_2 Y_4  , \\
   \alpha_1 &= X_1^3 + e^{(1)} X_1^2 X_2 ,   \\ 
   \alpha_2 &=  g^{(2)} X_1 X_2^2 + X_2^3 , \\
   \beta_1 &= Y_3^3 + h^{(1)} Y_3^2 Y_4 + l^{(1)} Y_3 Y_4^2 ,\\
   \beta_2 &=   h^{(2)} Y_3^2 Y_4 + l^{(2)} Y_3 Y_4^2 + Y_4^3 
\end{align*}
are homogeneous polynomials with coefficients in $\mathbb{C}$.  
We have a natural inclusion 
$\mathbb{C}^{14} = \{ (a^{(1)}, b^{(1)}, \cdots , l^{(2)}) \} 
    \hookrightarrow \mathbb{C}^{26}=\{ (a_j^{(i)})\} $, 
since the $F_i$'s above are special cases of $\tilde{F_i}$'s. 
We put $B_0 = \mathbb{C}^{14} \cap \tilde{B_0}$, 
and denote by $\psi : \mathcal{Y}_0 \to B_0$ and 
$\varphi : \mathcal{X}_0 \to B_0$, the pull-back of 
$\tilde{\mathcal{Y}} \to \tilde{B}$ and 
that of $\tilde{\mathcal{X}} \to \tilde{B}$, respectively. 
Now we are ready to prove Theorem \ref{thm:moduli}.
We use the same method as in \cite[Theorem 2. 11]{pg1k2=2} 
and \cite[Theorem 2.3]{globalpg1k2=1}.

\textbf{Proof of Theorem \ref{thm:moduli}.} 
Let $\mathcal{M}$ be the coarse moduli space for 
surfaces $X$'s as in Theorem \ref{thm:description}. 
By Theorem \ref{thm:description}, 
any surface $X$ as in  Theorem \ref{thm:description}  
corresponds to a fiber of the family $\tilde{X} \to \tilde{B}$, 
where $\tilde{B}$ is a non-empty Zariski open subset in $\mathbb{C}^{26}$. 
Thus by Tjurina's results on resolution of singularities 
(\cite{tjurinaresol}), all surfaces $X$'s as in  Theorem \ref{thm:moduli} 
are deformation equivalent, and their moduli space  $\mathcal{M}$ is 
irreducible. 
Meanwhile by the universality of the coarse moduli space, 
we have a natural morphism $B_0 \to \mathcal{M}$ corresponding 
to the family $\varphi: \mathcal{X}_0 \to B_0$. 
This morphism is dominant by Lemma \ref{lm:coefreduction} 
and its succeeding argument. 
By this together with Theorem \ref{thm:kuranishisp}, we obtain the 
unirationality of $\mathcal{M}$, since $B_0$ is a Zariski open subset 
in $\mathbb{C}^{14}$.  \qed

\section{The infinitesimal period map}   \label{scn:inftorelli}

In this section, we will give a proof for Theorem \ref{thm:infperiodmap}. 
Let $\psi : \mathcal{Y}_0 \to B_0$ and $\varphi : \mathcal{X}_0 \to B_0$ 
be the two analytic families given  
in Section \ref{scn:unirat}. 
For each 
$t = (a^{(1)}, b^{(1)}, \cdots , l^{(2)}) \in B_{0}$, 
the fibers $X= \varphi^{-1} (t)$ and $Y = \psi^{-1} (t)$ 
are a surface with 
invariants as in Theorem \ref{thm:moduli} and its universal cover,  
respectively. Note that $Y$ is a complete intersection in 
$W =\mathbb{P}^4$ defined by $F_{1}$ and $F_{2}$ as in 
(\ref{eqn:eqnunivcov}). 
We denote by $\pi : Y \to X $ and $\iota : Y \to W$ 
the natural projection and the natural inclusion, respectively.

Let $T_{B_0, t}$ be the holomorphic tangent space at 
$t \in B_0$. We denote by 
$\rho : T_{B_0, t} \to H^1 (\varTheta_X)$ and 
$\rho^{\prime} : T_{B_0, t} \to H^1 (\varTheta_Y)$ 
the Kodaira-Spencer map of $\varphi$ and that of $\psi$, 
respectively. 
In order to prove Theorem \ref{thm:infperiodmap}, 
it only suffices, by Theorem \ref{thm:kuranishisp} 
and the equality $\dim T_{B_0, t} = 14$, to show the injectivity of 
$\mu \circ \rho$ for general $t \in B_0$, 
where $\mu$ is the morphism given in Theorem \ref{thm:infperiodmap}. 
Note that the composite 
$\mu \circ \rho 
$ 
corresponds to 
the infinitesimal period map of $\varphi$. 
Let $\omega \in H^0 (\varOmega_X^2)$ be a non-zero 
holomorphic $2$-form on $X$ such that $\pi^* \omega$ 
corresponds to 
the section $W_0 \in H^0 (\varOmega_Y^2)$ 
in Remark \ref{rem:galoisgp}. 
Since $p_g (X) =1$, the kernel of $\mu \circ \rho $ is 
equal to that of the morphism 
$T_{B_0, t} \ni \xi \mapsto  
((\mu \circ \rho)(\xi))(\omega) 
\in H^1(\varOmega_X^1)$. 
Meanwhile we have the following commutative diagram:
\[ 
\begin{CD}  
T_{B_0, t}   @>\text{$\rho$}>>  H^1(\varTheta_X)  @>\text{$\omega \times$}>>
H^1(\varOmega^2_X \otimes \varTheta_X) \simeq H^1 (\varOmega^1_X)\\
@VV\text{$\parallel \mathrm{id}$}V  @VV\text{$$}V @VV\text{$$}V\\
T_{B_0, t}   @>\text{$\rho^{\prime}$}>>  H^1(\varTheta_Y)  
@>\text{$W_0 \times$}>> H^1(\varOmega^2_Y \otimes \varTheta_Y) 
\simeq H^1 (\varOmega^1_Y) , 
\end{CD}
\]
where the vertical morphisms 
$H^1(\varTheta_X) \to H^1(\varTheta_Y)$ and 
$H^1(\varOmega^2_X \otimes \varTheta_X) 
\to H^1(\varOmega^2_Y \otimes \varTheta_Y)$ 
are natural inclusions induced by the decomposition
of $\pi_* \varTheta_Y$ and 
$\pi_* (\varOmega_Y^2 \otimes \varTheta_Y)
\simeq \pi_* \varOmega^1_Y$ 
associated with the action of the Galois group $G=\mathrm{Gal} (Y/X)$.
Note that 
$((\mu \circ \rho ) (\xi))(\omega) = 
((\omega \times) \circ \rho ) (\xi)$ 
for any $\xi \in T_{B_0 , t}$. 
Thus in order to prove the injectivity of $\mu \circ \rho$, 
we only need to show that of the morphism 
$(W_0 \times) \circ \rho^{\prime} : 
T_{B_0, t} \to H^1(\varOmega^2_Y \otimes \varTheta_Y) $.  

Let us prove the injectivity of 
$(W_0 \times) \circ \rho^{\prime}$ for general $t \in B_0$. 
We denote by $R \simeq \oplus^{\infty}_{n=0} R_n$
and $R_n$ the graded ring 
$\mathbb{C} [W_0, X_1, X_2, Y_3, Y_4] / \langle F_1, F_2 \rangle$ 
and its homogeneous part of degree $n$, respectively.
This graded ring $R$ is naturally isomorphic to the canonical ring 
of $Y$. For each $m=0, 1, -1$, we denote by $R_n^{(m)}$ the 
set of all $F \in R_n$ satisfying 
\[
F(W_0,\varepsilon X_1,     \varepsilon X_2, 
       \varepsilon^{-1}Y_3, \varepsilon^{-1} Y_4)
=\varepsilon^m F(W_0,X_1,X_2,Y_3,Y_4) .
\]
This space $R_n^{(m)}$ corresponds to the 
eigenspace $H^0(\mathcal{O}_X (nK_X - mT_0))$ 
via the isomorphism (\ref{eql:decomposition}).

We have a natural exact sequence 
$
0 \to \varTheta_Y \to \iota^* \varTheta_W \to 
\mathcal{O}_Y (3)^{\oplus 2} \to 0
$
of $\mathcal{O}_Y$-modules.
By the similar argument as in Catanese \cite{pg1k2=1} 
and Oliverio \cite{torellik2=2}, we obtain, from this short exact sequence, 
the following commutative diagram: 
\begin{equation} \label{dgm:longexact}
\begin{CD}  
   R_1^{\oplus 5}    @>\text{$$}>> R_3^{\oplus 2} @>\text{$$}>> 
   H^1(\varTheta_Y) @>\text{$$}>> 0 \\
   @VV\text{$W_0 \times$}V  @VV\text{$W_0 \times$}V @VV\text{$W_0 \times$}V
   @VV\text{$$}V  \\
   R_2^{\oplus 5}    @>\text{$\delta$}>> R_4^{\oplus 2} @>\text{$$}>> 
   H^1(\varOmega^2_Y \otimes \varTheta_Y) @>\text{$$}>> \mathbb{C}
   @>\text{$$}>> 0 , 
\end{CD}
\end{equation}
where both of the horizontal sequences are exact, and the 
morphisms $R_1^{\oplus 5} \to R_3^{\oplus 2}$ and 
$\delta : R_2^{\oplus 5} \to R_4^{\oplus 2}$ 
are given by the matrix 
\[
\begin{pmatrix}
       \frac{\partial F_1}{\partial W_0} & 
       \frac{\partial F_1}{\partial X_1} &
       \frac{\partial F_1}{\partial X_2} &
       \frac{\partial F_1}{\partial Y_3} &
       \frac{\partial F_1}{\partial Y_4} \\
          \frac{\partial F_2}{\partial W_0} &
          \frac{\partial F_2}{\partial X_1} &
          \frac{\partial F_2}{\partial X_2} &
          \frac{\partial F_2}{\partial Y_3} &
          \frac{\partial F_2}{\partial Y_4}  
\end{pmatrix}
.
\]

Let $A^{\prime} : T_{B_0 , t} \to R_3^{\oplus 2}$ be the morphism 
given by 
$\frac{\partial}{\partial t} \mapsto 
     (\frac{\partial }{\partial t}F_1,  \frac{\partial }{\partial t}F_2)$, 
that is, the morphism giving the infinitesimal displacement of 
the deformation $\psi : \mathcal{Y}_0 \to B_0$ of 
the submanifold $Y \subset W$.  
Since the composite $(W_0 \times) \circ A^{\prime}$ maps 
$T_{B_0, t}$ into the subspace 
$R_4^{(0) \oplus 2} \subset R_4^{\oplus 2}$, 
we obtain a restriction 
$A: T_{B_0, t} \to R_4^{(0) \oplus 2}$ 
of $(W_0 \times) \circ  A^{\prime}$. 
We put $V = R_2^{(0)} \oplus R_2^{(1)} \oplus R_2^{(1)} 
\oplus R_2^{(-1)} \oplus R_2^{(-1)} \subset R_2^{\oplus 5}$,  
and denote by $C : V \to R_4^{(0) \oplus 2}$ 
%
%
the restriction of $\delta : R_2^{\oplus 5} \to  R_4^{\oplus 2}$
to this subspace.
Then from the commutative diagram (\ref{dgm:longexact}), 
we infer the equality 
$\ker ((W_0 \times) \circ \rho^{\prime}) = A^{-1} (C(V))$, 
where $C(V)$ is the image of the morphism $C$.

Let $M^{\prime}$ be a   
$26$-dimensional 
subspace of $R_4^{(0) \oplus 2}$ 
spanned by the following linearly independent elements: 
%
%
\begin{align} \label{eql:mprimebase}
&(W_0^4, 0), \quad (W_0 X_1 X_2^2, 0), \quad
(W_0 Y_3^3, 0), \quad (W_0 Y_4^3, 0), \quad 
(X_1^2 Y_3^2, 0),&  \notag \\
&(X_1^2 Y_3 Y_4, 0), \quad (X_1^2 Y_4^2, 0), \quad 
(X_1 X_2 Y_3^2, 0), \quad (X_1 X_2 Y_3 Y_4, 0),&  \notag \\
&(X_1 X_2 Y_4^2, 0), \quad 
(X_2^2 Y_3^2, 0), \quad (X_2^2 Y_3 Y_4, 0), \quad (X_2^2 Y_4^2, 0),&  \notag \\
&(0, W_0^4), \quad (0, W_0 X_1^2 X_2), \quad
(0, W_0 Y_3^3), \quad (0, W_0 Y_4^3), \quad 
(0, X_1^2 Y_3^2),& \notag \\
&(0, X_1^2 Y_3 Y_4), \quad (0, X_1^2 Y_4^2), \quad 
(0, X_1 X_2 Y_3^2), \quad (0, X_1 X_2 Y_3 Y_4),& \notag \\
&(0, X_1 X_2 Y_4^2), \quad 
(0, X_2^2 Y_3^2), \quad (0, X_2^2 Y_3 Y_4), \quad (0, X_2^2 Y_4^2). &  
\end{align}
Then, denoting the image of $A : T_{B_0, t} \to R_4^{(0) \oplus 2}$ 
by $M$, we have $R_4^{(0) \oplus 2} = M \oplus M^{\prime}$. 
Thus there exist two morphisms $D : V \to M$ and 
$D^{\prime} : V \to M^{\prime}$ such that $C = D + D^{\prime}$. 
Note that $C(V) \cap M \simeq D (\ker D^{\prime})$. 
By this together with the injectivity of $A$, we obtain 
\[
 \ker ((W_0 \times) \circ \rho^{\prime}) =  A^{-1} (C(V))
                                         \simeq D (\ker D^{\prime}). 
\] 
Meanwhile we have 
$\dim V = 25$ and 
\[
 (W_0^2, W_0 X_1, W_0 X_2, W_0 Y_3, W_0 Y_4) \in \ker C 
  = \ker D \cap \ker D^{\prime}.  
\]
Thus in order to prove the injectivity of 
$ (W_0 \times) \circ \rho^{\prime} : T_{B_0 , t} 
\to H^1(\varOmega_Y^1 \otimes  \varTheta_Y)$, 
we only need to show the equality $\mathrm{rank} D^{\prime} = 24$.  

So, in what follows, we will show $\mathrm{rank} D^{\prime} = 24$ 
for general $t \in B_0$. 
We employ the following base of $V$: 
%
%
%
\begin{align} \label{eql:vbase}
&(W_0^2)_1, \quad 
(X_1 Y_3)_1, \quad 
(X_1 Y_4)_1, \quad 
(X_2 Y_3)_1, \quad
(X_2 Y_4)_1,& \notag \\
&(W_0 X_1)_2, \quad 
(W_0 X_2)_2, \quad
(Y_3^2)_2, \quad 
(Y_3 Y_4)_2, \quad
(Y_4^2)_2,& \notag \\
&(W_0 X_1)_3, \quad
(W_0 X_2)_3, \quad
(Y_3^2)_3, \quad
(Y_3 Y_4)_3, \quad 
(Y_4^2)_3,& \notag \\
&(W_0 Y_3)_4, \quad
(W_0 Y_4)_4, \quad
(X_1^2)_4, \quad
(X_1 X_2)_4, \quad
(X_2^2)_4,& \notag \\
&(W_0 Y_3)_5, \quad
(W_0 Y_4)_5, \quad
(X_1^2)_5, \quad
(X_1 X_2)_5, \quad
(X_2^2)_5,&  
\end{align}
where, for each $u \in R_2$ and $1 \leq i \leq 5$, 
we denote by $(u)_i$ the element 
$(v_1, v_2, v_3, v_4, v_5) \in R_2^{\oplus 5}$ given by  
$v_i = u$, $v_j = 0$ $(j \neq i)$. 
Let $L_1$ be the $26 \times 25$ matrix of $D^{\prime}$ 
corresponding to the bases (\ref{eql:vbase}) of $V$
and (\ref{eql:mprimebase}) of $M^{\prime}$: i.e.  
\[
L_1 = 
\begin{pmatrix}
L_{1,1}^{(1)}& L_{1,2}^{(1)}& L_{1,3}^{(1)} \\
L_{2,1}^{(1)}& L_{2,2}^{(1)}& L_{2,3}^{(1)} \\
\end{pmatrix}, 
\]
where $13 \times 5$ matrixes $L_{1,1}^{(1)}$, $L_{2,1}^{(1)}$, 
and $13 \times 10$ matrixes $L_{1,2}^{(1)}$, $L_{2,2}^{(1)}$, 
$L_{1,3}^{(1)}$, $L_{2,3}^{(1)}$ are given by  
\begin{small}
\[ 
L_{1,1}^{(1)} = 
\begin{bmatrix}
 3& & & & \\
 0& & & & \\
 0& & & & \\
 0& & & & \\
 & a^{(1)}& & & \\
 & b^{(1)}& a^{(1)}& & \\
 & & b^{(1)}& & \\
 & c^{(1)}& & a^{(1)}& \\
 & d^{(1)}& c^{(1)}& b^{(1)}& a^{(1)}\\
 & & d^{(1)}& & b^{(1)}\\
 & & & c^{(1)}& \\
 & & & d^{(1)}& c^{(1)}\\
 & & & & d^{(1)}\\
\end{bmatrix}
, \quad
L_{2,1}^{(1)} =
\begin{bmatrix}
 3& & & & \\
 0& & & & \\
 0& & & & \\
 0& & & & \\
 & a^{(2)}& & & \\
 & b^{(2)}& a^{(2)}& & \\
 & & b^{(2)}& & \\
 & c^{(2)}& & a^{(2)}& \\
 & d^{(2)}& c^{(2)}& b^{(2)}& a^{(2)}\\
 & & d^{(2)}& & b^{(2)}\\
 & & & c^{(2)}& \\
 & & & d^{(2)}& c^{(2)}\\
 & & & & d^{(2)}\\
\end{bmatrix}, 
\] 
\end{small}
\begin{small}
\[
 L_{1,2}^{(1)} = 
\begin{bmatrix}
 -3& & & & & -e^{(1)}& & & & \\
 & 2e^{(1)}& & & & & & & & \\
 -3& & a^{(1)}& & & -e^{(1)}& 0& c^{(1)}& & \\
 & & & & b^{(1)}& & & & & d^{(1)}\\
 & & 3& & & & & e^{(1)}& & \\
 & & & 3& & & & & e^{(1)}& \\
 & & & & 3& & & & & e^{(1)}\\
 & & 2e^{(1)}& & & & & & & \\
 & & & 2e^{(1)}& & & & & & \\
 & & & & 2e^{(1)}& & & & & \\
 & & & & 0& & & & & \\
 & & & & 0& & & & & \\
 & & & & 0& & & & & \\
\end{bmatrix}, 
\]
\end{small}
\begin{small}
\[
 L_{2,2}^{(1)} = 
\begin{bmatrix}
 & -g^{(2)}& & & & & -3& & & \\
 & & & & & 2g^{(2)}& & & & \\
 & & a^{(2)}& & & & & c^{(2)}& & \\
 0& -g^{(2)}& & & b^{(2)}& & -3& & & d^{(2)}\\
 & & & & & & & 0& & \\
 & & & & & & & 0& & \\
 & & & & & & & 0& & \\
 & & & & & & & 2g^{(2)}& & \\
 & & & & & & & & 2g^{(2)}& \\
 & & & & & & & & & 2g^{(2)}\\
 & & g^{(2)}& & & & & 3& & \\
 & & & g^{(2)}& & & & & 3& \\
 & & & & g^{(2)}& & & & & 3\\
\end{bmatrix}, 
\]
\end{small}
\begin{small}
\[
L_{1,3}^{(1)} = 
\begin{bmatrix}
 & & -a^{(1)}& & -c^{(1)}& & & -b^{(1)}& & -d^{(1)}\\
 & & & c^{(1)}& * & & & & d^{(1)}& * \\
 3& & -a^{(1)}& & & h^{(1)}& 0& -b^{(1)}& & \\
 & l^{(1)}& & & -c^{(1)}& & & & & -d^{(1)}\\
 & & 3& & & & & h^{(1)}& & \\
 & & 2h^{(1)}& & & & & 2l^{(1)}& & \\
 & & l^{(1)}& & & & & & & \\
 & & & 3& & & & & h^{(1)}& \\
 & & & 2h^{(1)}& & & & & 2l^{(1)}& \\
 & & & l^{(1)}& & & & & & \\
 & & & & 3& & & & & h^{(1)}\\
 & & & & 2h^{(1)}& & & & & 2l^{(1)}\\
 & & & & l^{(1)}& & & & & \\
\end{bmatrix}, 
\]
\end{small}

\begin{small}
\[
 L_{2,3}^{(1)} =
\begin{bmatrix}
 & & -a^{(2)}& & -c^{(2)}& & & -b^{(2)}& & -d^{(2)}\\
 & & *& a^{(2)}& & & & *& b^{(2)}& \\
 & & -a^{(2)}& & & h^{(2)}& & -b^{(2)}& & \\
 0& l^{(2)}& & & -c^{(2)}& & 3& & & -d^{(2)}\\
 & & & & & & & h^{(2)}& & \\
 & & 2h^{(2)}& & & & & 2l^{(2)}& & \\
 & & l^{(2)}& & & & & 3& & \\
 & & & & & & & & h^{(2)}& \\
 & & & 2h^{(2)}& & & & & 2l^{(2)}& \\
 & & & l^{(2)}& & & & & 3& \\
 & & & & & & & & & h^{(2)}\\
 & & & & 2h^{(2)}& & & & & 2l^{(2)}\\
 & & & & l^{(2)}& & & & & 3\\
\end{bmatrix}.
\]
\end{small}
Here empty entries are zero. 
For general $t \in B_0$, 
we strike off the rows and the columns of $L_1$ meeting 
the following entries by doing the operations in this order: 
$(3,16)$,$(1,6)$,$(17,22)$,$(14,12)$,
$(2,7)$,$(4,17)$,$(15,11)$,$(16,21)$. 
Then we see $\mathrm{rank} L_1 = \mathrm{rank} L_2 +8$, 
where $L_2$ is the $18 \times 16$ matrix obtained by 
removing from $L_1$ the following: 
i) the rows and columns meeting the $8$ entries given above, 
and ii) the first column. 
Thus we only need to show $\mathrm{rank} L_2 = 16$ for general $t \in B_0$. 

Let $L_3$ be the $18 \times 16$ matrix obtained by 
specializing $L_2$ by $e^{(1)} = g^{(2)} = 0$. 
We can strike off the rows and the columns of $L_3$ meeting 
the following entries: 
$(1,5)$,$(2,6)$,$(3,7)$,$(16,8)$,$(17,9)$,$(18,10)$. 
Thus we see 
$\mathrm{rank} L_2 \geq \mathrm{rank} L_3 
= \mathrm{rank} L_4 + 6$ for general $t \in B_0$, 
where $L_4$ is the $12 \times 10$ matrix obtained 
by removing from $L_3$ the rows and columns meeting 
the $6$ entries above.    
Hence we only need to show $\mathrm{rank} L_4 =10$ 
for general $t \in B_0$. 

It now suffices to show $\det L_5 \neq 0$ for general $t \in B_0$, 
where the $10 \times 10$ matrix 
\begin{small}
\begin{equation} \label{eql:matnoteno6}
L_5 = 
\begin{bmatrix}
 c^{(1)}& & a^{(1)}& & & 3& & & h^{(1)}& \\
 d^{(1)}& c^{(1)}& b^{(1)}& a^{(1)}& & 2h^{(1)}& & & 2l^{(1)}& \\ 
 & d^{(1)}& & b^{(1)}& & l^{(1)}& & & & \\ 
 & & c^{(1)}& & & & 3& & & h^{(1)}\\ 
 & & d^{(1)}& c^{(1)}& & & 2h^{(1)}& & & 2l^{(1)}\\ 
 b^{(2)}& a^{(2)}& & & 2h^{(2)}& & & 2l^{(2)}& & \\
 & b^{(2)}& & & l^{(2)}& & & 3& & \\ 
 c^{(2)}& & a^{(2)}& & & & & & h^{(2)}& \\ 
 d^{(2)}& c^{(2)}& b^{(2)}& a^{(2)}& & 2h^{(2)}& & & 2l^{(2)}& \\ 
 & d^{(2)}& & b^{(2)}& & l^{(2)}& & & 3& \\      
\end{bmatrix}
\end{equation}
\end{small}
is the one obtained by 
removing from $L_4$ its $6$-th and $7$-th rows.  
But, when we compute $\det L_5$ by the definition of 
the determinant, the monomial 
$(d^{(2)})^2 (a^{(1)})^2 (l^{(2)})^2  h^{(2)} (h^{(1)})^2 l^{(1)}$ 
appears only once, i.e., from the term passing the entries 
$(9,1)$, $(10,2)$, $(1,3)$, $(2,4)$, $(7,5)$, 
$(3,6)$, $(5,7)$, $(6,8)$, $(8,9)$, and $(4,10)$ of $L_5$. 
Thus, for general $t \in B_0$, we have $\det L_5 \neq 0$, 
and hence $\mathrm{rank} D^{\prime} = 24$,  
which completes the proof of Theorem \ref{thm:infperiodmap}.

\begin{flushright}
\begin{minipage}{24em}
Masaaki Murakami \\
Department of Mathematics, 
Faculty of Science, \\
Kyoto University,
Kyoto 606-8502, Japan \\
\end{minipage}
\end{flushright}
 

\end{document}